# The oldest example of $\pi \approx 3\frac{1}{8}$ in Sumer: Calculation of the area of a circular plot

Kazuo MUROI

§1. Two formulas

The circle is one of the fundamental figures in elementary geometry and so it is natural that the ancients would be interested in calculating its area. It is well known that the Babylonians of the Old Babylonian period (ca. 2000-1600 BCE) used two formulas for the area enclosed by a circle whose circumference is c:

$$\frac{1}{4\pi}c^2 \approx \frac{1}{12}c^2 = \mathbf{0;5}\ c^2\ (\pi \approx 3)\text{ and}$$

$$\frac{1}{4\pi}c^2 \approx \frac{2}{25}c^2 = \mathbf{0;4,48}\ c^2\ (\pi \approx 3\frac{1}{8})$$

The two constants 0;5 and 0;4,48 are also listed in some tables of constants of those days[1]:

5 ša gán-gúr "0;5 is that of a circular field."

4,48 ša ku-bu-ur-re-e giš "0;4,48 is that of the thickness of a log."

4,48 é-kišib še "0;4,48 is (that of) the cylindrical storehouse of barley."

In mathematical problems the first formula frequently occurs but the second has only been found once so far, in a problem which concerns the volume of a cylindrical log.[2]

A few comments should be made about the final item in the list and particularly the Sumerian word é-kišib-(ba), which was used from the Sargonic period (ca. 2340-2200 BCE) onwards. The literal meaning is "house of a cylinder seal ", and it has occasionally been mistranslated as "sealed storehouse". As it is itemized in a list of mathematical constants, this suggests that the word has something to do with a circular plot. Therefore, I have translated it as "cylindrical storehouse" (see above).

§2. Circular plot in a field plan

In 1903 details of a fragmentary tablet which was part of a field plan of a city were published by F. Thureau-Dangin.[3]  The tablet was excavated from the ruins of Sumerian city Girsu and came from the Sargonic period. Although it appears that no one has shown any interest in it since Thureau-Dangin's publication, the tablet may be of great importance from the viewpoint of the history of mathematics, because a circular plot drawn in the tablet clearly shows the use of π ≈ 3 + 1/8 (= 3;7,30).  Inside the circle the area is written down as follows:

   7 sar 10 gín, where 1 sar = 60 gín = 1 nindan² ≈ 36 m².   See Fig.1.

If the diameter of the circle is 3 nindan 1 šu-dù-a (= 3;1,40 nindan) ,where 1 nindan = 36 šu-dù-a ≈ 6 m, and π ≈ 3;7,30, the circumference and its square would be:

   c = 3;7,30 · 3;1,40 = 9;27,42,30

   c² = (9;27,42,30)² = 1,29;31,32,45,6,15.

Therefore, the area of the circular plot is:

$$\frac{1}{4\pi}c^2 \approx 0;4,48 \cdot 1,29;31,32,45,6,15$$

   = 7;9,43,25,12,30

   ≈ 7;10 (sar)

= 7 (sar) 10 (gín).

It is probable that 3 nindan (≈ 18 m) is the inside diameter of the base of a cylindrical storehouse, é-kišib-(ba), and 1/2 šu-dù-a (≈ 8.3 cm) is the width of the surrounding wall or fence.   See Fig. 2.   The length unit šu-dù-a is particular to this period   or Pre-Sargonic period (ca. 2600-2340 BCE) and it was never used in later periods.

§3. Origin of the two formulas

Using the method of exhaustion, an elementary approach to the integral calculus, and the method of a proof by contradiction, Archimedes (287-212 BCE) demonstrated the theorem that the area enclosed by a circle of radius r is equal to that of a triangle whose base and height are equal to the circumference and the radius respectively.[4]   It is conceivable that the Sumerians had already known the theorem above without proving it precisely. If so, since $r = c/2\pi$, they needed to choose a suitable number for an approximate value of $\pi$ whose reciprocal should be finite in their sexagesimal place value notation.   Out of many reciprocal tables such as:

$(2^n, 2^{-n})$, $(2^n \cdot 3, 2^{-n} \cdot 3^{-1})$, $(2^n \cdot 3^2, 2^{-n} \cdot 3^{-2})$, ---, $(2^n \cdot 5, 2^{-n} \cdot 5^{-1})$, $(2^n \cdot 5^2, 2^{-n} \cdot 5^{-2})$, ---

$(3^n, 3^{-n})$ --- and so on for n = 1   2   3 ---- 30,[5]

they chose the two handy pairs $(3, 3^{-1})$ and $(2^7 \cdot 3^2, 2^{-7} \cdot 3^{-2})$ or (19,12, 0;0,**3,7,30**), yielding $1/4\pi \approx 0;5$ and $1/4\pi \approx 0;4,48$ respectively. Moreover one pair, 3,8,22,48,24,10 (= $2 \cdot 5^{13}$) and 0;0,0,0,0,0,19,6,37,4,16,53,45,36 (= $2^{-1} \cdot 5^{-13}$) may give a more accurate approximation to $\pi$ because :

$$3;8,22,48,24,10 \approx 3.13967,$$

but it is obviously very hard to use it for practical calculation. It seems that there were only two choices of reciprocal pairs available to the Sumerians to get an approximation to $\pi$.

As I have already pointed out elsewhere,[6] it is possible that the Babylonians used quite an accurate approximation to $\pi$, that is $\pi \approx 3;9$ (=3.15), for the area of a segment, but this value is not fit for $c^2/4\pi$ because its reciprocal cannot be obtained by their sexagesimal notation. Thus the choice of the number for $\pi$ in the formula $c^2/4\pi$ may have depended heavily on whether its reciprocal is finite or not in sexagesimal number system.

§4. Conclusion

We have confirmed that the formula $0;4,48\ c^2$ for the area of a circle, which occurs in a few Babylonian mathematical texts, can be traced back to at the latest the twenty-third century BCE . This fact, together with some others,[7] leads us to the conclusion that the basis of Babylonian mathematics had been established by the end of the third millennium BCE. Further studies of the cuneiform tablets of Pre-Sargonic period and Sargonic period may give us historically important facts concerning the origin of Babylonian mathematics, since, in my judgment, some of the tablets have not been studied enough from the mathematical point of view. If we advance our studies along these lines, we will be able to describe "Sumerian mathematics" in more detail someday.

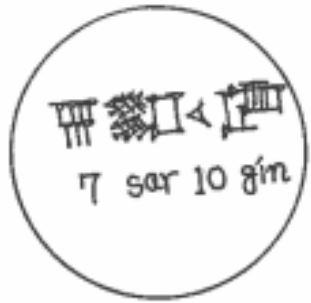 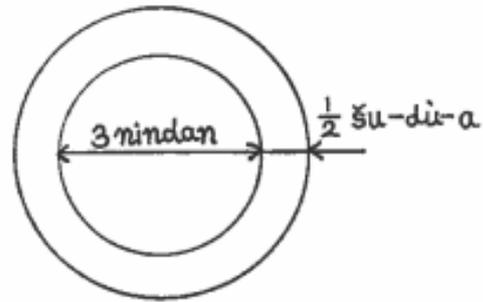

Fig.1  Fig.2

Notes

(1) O. Neugebauer and A. Sachs, *Mathematical Cuneiform Texts* (= MCT), 1945, pp. 132-136, YBC 5022 lines 20, 58, and 61.

(2) MCT, pp. 57-59, YBC 8600.

(3) F. Thureau-Dangin, *Recueil de tablettes chaldéennes*, 1903, no. 150.

(4) M. R. Cohen and I. E. Drabkin, *A Source Book in Greek Science*, 1948, pp. 59-61.

(5) K. Muroi, How to Construct a Large Table of Reciprocals of Babylonian Mathematics, arXiv: 1401. 0065 [moth. HO]

(6) K. Muroi, Mathematics Hidden Behind the Practical Formulae of Babylonian Geometry, in G. J. Selz, ed., *The Empirical Dimension of Ancient Near Eastern Studies*, 2011, pp. 149-157.

(7) For example, see my paper:

K. Muroi, The oldest example of compound interest in Sumer: Seventh power of four-thirds, arXiv: 1510. 00330 [math. HO].